\documentclass{article}

\usepackage{PRIMEarxiv}

\usepackage[utf8]{inputenc} 
\usepackage[T1]{fontenc}    
\usepackage{hyperref}       
\usepackage{url}            
\usepackage{booktabs}       
\usepackage{amsfonts}       
\usepackage{nicefrac}       
\usepackage{microtype}      
\usepackage{lipsum}
\usepackage{fancyhdr}       
\usepackage{graphicx}       
\graphicspath{{media/}}     
\usepackage{amsmath}

\usepackage{algorithm, setspace}
\usepackage{algpseudocode}
\usepackage{cleveref}
\usepackage{bm}
\usepackage{gensymb}
\usepackage{adjustbox}

\newlength\myindent
\title{Design and Guidance of a Multi-Active Debris Removal Mission}

\author{
  Minduli C. Wijayatunga , Roberto Armellin , Harry Holt, Laura Pirovano  \\
  The University of Auckland,   \\
  Auckland 1010,  \\
  New Zealand\\
  \texttt{mwij516@aucklanduni.ac.nz} \\
   \And
  Aleksander~A.~Lidtke \\
   Astroscale Japan Inc., \\
  1-16-4-16 Kinshi, Sumida-ku,   \\
  Tokyo, Japan\\
  \texttt{a.lidtke@astroscale.com} \\
}

\begin{document}
\maketitle

\begin{abstract}
Space debris have been becoming exceedingly dangerous over the years as the number of objects in orbit continues to rise. 
Active debris removal (ADR) missions have garnered significant attention as an effective way to mitigate this collision risk. This research focuses on developing a multi-ADR mission that utilizes controlled reentry and deorbiting. The mission comprises two spacecraft: a Servicer that brings debris down to a low altitude and a Shepherd that rendezvous with the debris to later perform a controlled reentry. A preliminary mission design tool (PMDT) is developed to obtain time or fuel optimal trajectories for the proposed mission while taking the effect of $J_2$, drag, eclipses, and duty ratio into account. The PMDT can perform such trajectory optimizations within computational times that are under a minute. Three guidance schemes are also studied, taking the PMDT solution as a reference, to validate the design methodology and provide guidance solutions for this complex mission profile.
\end{abstract}


\section{Introduction}
The space environment in low Earth orbit (LEO) is increasingly populated with space debris. As a result, the average rate of debris collisions has increased to four or five objects per year \cite{Maestrini2021}. Most debris are artificial objects, including derelict satellites, discarded rocket stages, and fragments originating from collisions. As satellites become increasingly essential to daily life, more and more satellites are added to expand space-enabled services. However, additional launches increase the risk of collision for all satellites as they further saturate space with objects, thereby endangering the critical space infrastructure. A collision in space can create debris that can collide with other space objects and generate more debris.
This cascading effect is known as the \lq\lq Kessler Syndrome\rq\rq, named after D.J. Kessler \cite{Kessler1, Kessler2}. Kessler et al. \cite{Kessler2} discussed the frequency of collisions and their consequences, describing standard mitigation techniques for the first time. Then, Pelton \cite{pelton} discussed the cascading effect of collisions and the international standards for debris mitigation and space traffic management. He also gave estimates for the number of orbital debris at the time to be around six metric tonnes in mass and 22000 in number. Several events in recent history have caused significant additions to the space debris population. These include the anti-satellite missile tests in 2007 and 2021, and the collision of Iridium 33 and Kosmos 2251 in 2009 \cite{surveydebris,potter_2021}.

\par 
Active debris removal (ADR) is the process of removing derelict objects from space, thus minimizing the build-up of unnecessary objects and lowering the probability of on-orbit collisions that could fuel the \lq\lq collision cascade \rq\rq  \cite{bonnal2013active,liou2010parametric}. ADR has gained traction in the past two decades, leading to numerous studies and implementations of potential debris removal missions and technologies. The ELSA-d mission designed by Astroscale was launched in March 2021 and has successfully tested rendezvous algorithms needed for ADR and a magnetic capture mechanism needed to remove objects carrying a dedicated docking plate at the end of their missions  \footnote{\url{https://astroscale.com/astroscales-elsa-d-mission-successfully-completes-complex-rendezvous-operation/}}. The RemoveDebris mission by the University of Surrey is another project that demonstrated various debris removal methods, including harpoon and net capture \cite{Forshaw2020TheAS}. The CleanSpace-1 mission by the European Space Agency (ESA) aims to deorbit a~112~kg upper stage of the Vega rocket \footnote{\url{https://www.esa.int/Space_Safety/ClearSpace-1}}.\par
While individual removals are essential stepping stones towards ADR implementation, a deployment on a larger scale, targeting more objects, might be necessary \cite{liou2010parametric, white2014}. In order to make it financially feasible, each ADR Servicer might need to remove more than one object and use mass-efficient low-thrust electric propulsion (EP). This combination of long, EP-based transfers and complex vehicle paths rendezvousing with multiple moving targets presents a difficult optimization challenge that has to be addressed at the design stage of ADR missions. Due to the large number of potential ADR targets to be visited, transfers between consecutive mission orbits need to be analyzed quickly to enable design iteration and parametric studies. 

This paper introduces a novel multi-ADR removal mission concept that involves a two-spacecraft system. On request, the system is able to provide contact-based debris removal through a rendezvous and deorbit process. One spacecraft - called the Servicer - is reused for multiple debris, allowing the mission costs to remain low. The other spacecraft - the Shepherd - performs coupled reentry with the debris, so the reentry process can be controlled adequately, thus complying with the ground-casualty risk requirements. The majority of the mission utilizes electric propulsion. This paper is dedicated to discussing the proposed mission in detail and developing a mission design tool to simulate multi-ADR tours accurately and efficiently.  \par 

To this end, a preliminary mission design tool (PMDT) is developed to optimize both fuel consumption and time of flight of multi-target missions while taking the effect of $J_2$, eclipses, and duty ratio into account. PMDT utilizes $J_2$ to achieve RAAN changes, in order to reduce the fuel consumption of the mission. \par 
The PMDT extends the traditional Edelbaum method by introducing the contribution of drag and duty ratio. 
Then drift orbits are used for matching RAAN when required, as discussed in \cite{6}. Lastly, the altitude and inclination of the drift orbits are optimized to obtain either time or fuel optimal trajectories. The sequence of targets can also be treated as an optimization variable in the PMDT, however, it was treated as a constant for the examples given in this study. 
\par 
Our approach shares similarities with the Multidisciplinary desigN Electric Tug tool (MAGNETO) developed in \cite{Rimani} as well as the work by Viavattene et al. in \cite{Viavattene2022DesignOM}. However, it takes the presented models further by taking duty ratios into account and a more accurate description of eclipses and drag. Furthermore, the tool considers mission-specific constraints and uses an optimizer to perform rapid design iterations and parametric studies of the proposed multi-ADR mission. \par

Three guidance laws are introduced to assess the accuracy of the models adopted in the PMDT. Ruggiero et al. \cite{e7}  developed a series of closed-loop guidance laws based on the Gauss form of Lagrange Planetary laws. Locoche \cite{slim}  developed a guidance law based on Lyapunov feedback control known as the $\Delta v$-Law to supplement preliminary mission design tools. Finally, Petropoulos \cite{Petropoulos2005} developed one of the most versatile and well-known control laws - the Q-Law - which is also based on Lyapunov control. These three approaches are here used to optimally track the transfers computed by the PMDT, thus validating its key assumptions and providing a possible way to fly the missions. 
\par 
The remaining sections of the paper are organized as follows. Firstly, the mission concept of operations is presented. Then, the design of PMDT is discussed and is used to generate optimal debris removal trajectories at high computational speeds. Thirdly, guidance schemes implemented on the PMDT outcomes are discussed. The results section provides an example trajectory optimization solution for both a time and mass optimal multi-ADR mission. Lastly, the paper's outcomes are summarized, and conclusions are drawn regarding the method's usefulness. 

\section{Concept of Operations of the Multi-ADR Mission}
\begin{figure}[hbt!]
    \centering
    \includegraphics[width = \textwidth]{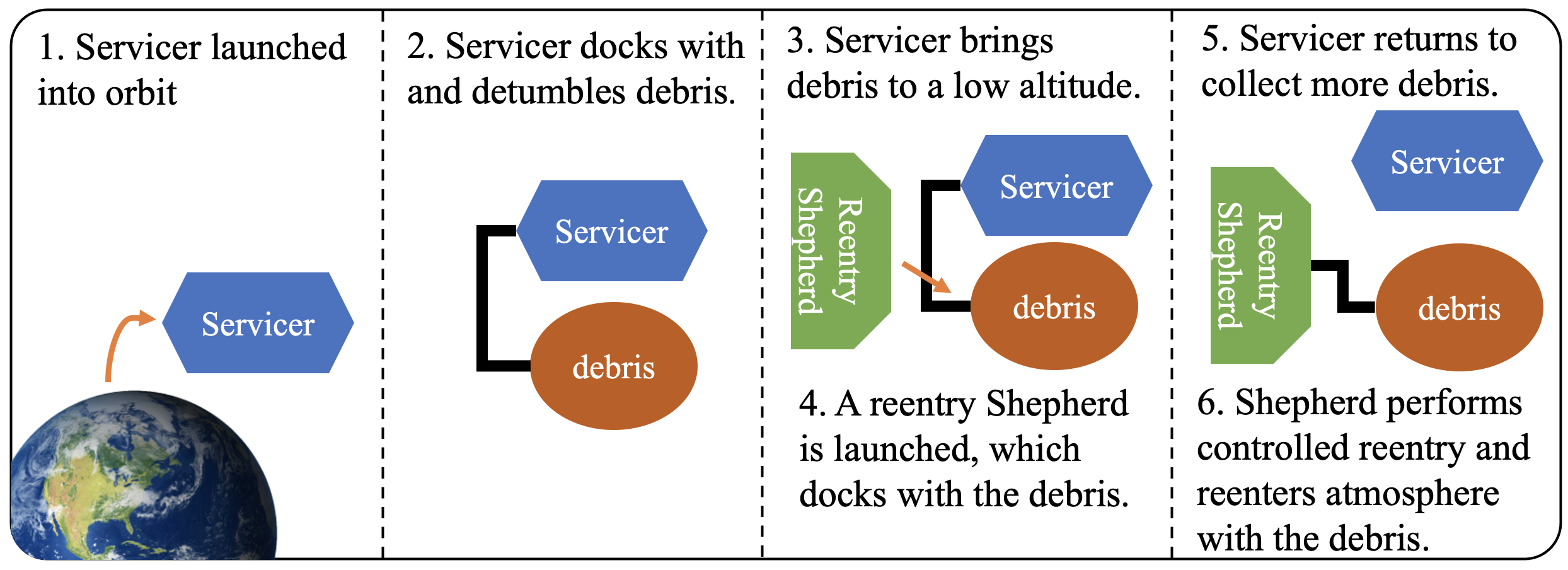}
    \caption{Mission architecture of the multi-ADR mission}
    \label{ADR}
\end{figure}

The proposed multi-ADR mission architecture is shown in \Cref{ADR}, where two spacecraft are involved in the debris removal process. A Servicer is used to approach and rendezvous with the debris. Once rendezvoused, the Servicer brings the object down to a low altitude orbit ($\approx 350$~km). The debris is then handed over to a Reentry Shepherd, which docks with the debris and performs a~controlled reentry on its behalf. Controlled reentry reduces the casualty risk posed by removing the debris, which is desirable because the ADR targets are, by definition, large and thus contain components likely to survive the reentry. The Servicer shall be reused for several debris removals, while each Reentry Shepherd can only be used once as it burns up while deorbiting the debris.

\par

The proposed mission architecture can perform multi-ADR services significantly cheaper than those that use coupled deorbiting and controlled reentry systems. When the deorbiting and reentry functionality are installed on a single spacecraft, it cannot be reused, which leads to higher mission costs.
Furthermore, depending on the debris features and the governing regulations, the requirements of controlled reentry for each debris will differ. The Reentry Shepherd allows the required flexibility for controlled reentry of large-scale objects. \par 
The separation of the two systems ensures that the Servicer is not obligated to carry re-entry-related hardware, thus saving mass for fuel. It also warrants that the Shepherd does not need to perform extensive orbital changes or have a long lifetime in space which, in turn, reduces its size and cost.

\par 
Trajectory optimization for this mission through traditional means is non-intuitive, as multiple transfer arcs and targets are involved. Hence, a trajectory optimization tool capable of developing suboptimal solutions with limited computational capacity is developed in the following section.

\section{Methodology}

\subsection{Design of the Preliminary Mission Design Tool}

\par 
In the 1960s, an analytical solution for the transfer between two inclined circular orbits under continuous thrust was developed by Theodore N. Edelbaum \cite{2}. While the transfer arcs developed were both time and fuel optimal, they were obtained under the assumption of continuous thrust and the lack of other perturbations such as $J_2$ and air resistance. Several studies were conducted following Edelbaum's work to include the effect of discontinuous thrust and orbital perturbations on the problem dynamics. Colasurdo and Casalino \cite{9} extended Edelbaum's analysis to compute optimal quasi-circular transfers while considering the effect of the Earth's shadow, and Kechichian \cite{3} developed a method for calculating coplanar orbit-raising maneuvers taking eclipses into account while constraining the eccentricity to zero. However, both \cite{9,3} could only provide suboptimal solutions, as they utilized thrust steering to maintain zero eccentricity. In 2011, Kluever \cite{4} further extended Kechichian's method into a semi-analytic method that considers the effect of $J_2$ and Earth-shadow arcs. This method used Edelbaum-based orbital elements to compute the Earth shadow arc during the transfer. However, it failed to consider the effect of air resistance, which is of crucial value for LEO transfers. In 2019, Cerf \cite{6} proposed utilizing $J_2$ to achieve right ascension of ascending node (RAAN) changes during transfers to reduce fuel consumption while keeping the time of flight constant and not taking the effect of eclipses and air resistance into account. The PMDT is developed to unite the ideas given in \cite{2,6,4} and take them a step further by considering air resistance and duty ratio. \par 
The PMDT first calculates the time of flight and fuel expenditure of a single transfer using Edelbaum's method described in \cite{2}. Then, additions to the classical Edelbaum method- creating the Extended Edelbaum method (Algorithm \ref{alg:1})- are made such that the effect of atmospheric drag, engine duty ratio, and solar eclipses are taken into account. Thirdly, a RAAN matching algorithm (Algorithm \ref{alg:3}) that does not utilize fuel to make RAAN adjustments is implemented to make transfers cheaper. This is achieved by introducing an intermediate drift orbit where the Servicer can utilize the effect of $J_2$ perturbations to reach the desired RAAN. Lastly, this process is introduced into an optimization scheme (Figure \ref{opt}) where the launch time and the drift orbits involved can be optimized to achieve the minimum time of flight or the minimum fuel expenditure. \par 



\subsubsection{Extended Edelbaum Method}
The Extended Edelbaum method is a version of the classical Edelbaum method adapted to take the effect of atmospheric drag, solar eclipses, and duty ratios into account. This method is detailed in \Cref{alg:1}. Note that the Extended Edelbaum method only ensures that a desired semi-major axis and inclination are reached. 
\par 

\begin{algorithm}[hbt!]
\caption{Extended Edelbaum Method}\label{alg:1}
\begin{algorithmic}
\Require Initial and final orbital velocity ($V_0,V_f$), change in inclination ($\Delta i$), thrust acceleration ($f$)
\State Calculate $\Delta v$ and mission time of flight ($t_f$) using Eq. \eqref{eq1} and \eqref{eq2}.
\begin{align} 
  \Delta v_{\text{total}} &= \sqrt{V_0^2 + V_f^2 - 2V_0 V_f \cos(\pi/2 \Delta i)}  \label{eq1}\\ 
    t_f &= \frac{ \Delta v_{\text{total}}}{f} \label{eq2}
\end{align}
\State Calculate the initial yaw steering angle $\beta$, defined in the plane normal to the orbit plane, using Eq. \eqref{eq3}. 
\begin{equation}\label{eq3}
\tan{\beta_0} = \frac{\sin{(\pi/2 \Delta i )}}{\frac{V_0}{V_f} - \cos(\pi/2 \Delta i )}
\end{equation}
\State Discretise the time of flight ($t_f$) into N segments and compute the semi-major axis, inclination and $\Delta v$ per segment using Eq. \eqref{eq4}, \eqref{eq5}, and \eqref{eq2}, respectively.
\begin{align} 
a(t) &= \frac{\mu}{V_0^2 + f^2t^2 - 2V_0ft\cos(\beta_0)}
 \label{eq4}\\ 
i(t) &= i_0 +\textrm{sgn}{(i_f-i_0)}\frac{2}{\pi}\left[\tan^{-1}\left(\frac{f t - V_0 \cos{\beta_0}}{V_0 \sin{\beta_0}}\right) + \frac{\pi}{2}-\beta_0 \right] \label{eq5}\end{align}
\For{$k = 1:N$}
    \State Calculate sunlit time during a single orbit ($w_{ecl}$), using eclipse time formulation in \cite{5}.
    \State Calculate the fraction of thrust time per orbit ($w$).
     \begin{equation}
    w =\min[DR, w_{ecl}] \ \text{where}  \ DR: \text{Duty Ratio}
    \end{equation}
    \State Compute the new transfer time using Eq. \eqref{eq6}. 
    \begin{equation}\label{eq6}
    t_{k+1} = t_k+ \frac{\Delta v_{k+1}- \Delta v_k }{f_k w_k}
    \end{equation}
    \State Calculate drag acceleration at $t_{k+1}$ and $t_k$ using   Eq. \eqref{eq7}. 
    \begin{equation}\label{eq7}
    a_{drag}=-\frac{1}{2} \frac{\rho C_d A v^2}{m}    
    \end{equation}
    \State $\rho, C_d, A, v$ and $m$ represent the air density, drag coefficient, frontal area, velocity, and mass. 
    \State Calculate $a_{k+1, drag}$ and $a_{k, drag}$ corresponding to each drag acceleration using Eq. \eqref{eq4}. 
    \If{$|a_{k+1, drag} - a_{k, drag}| > \Delta a$} 
    \State $a_{k+1} = a_{k+1} + |a_{k+1, drag} - a_{k, drag}| $
    \State Go back to the first step of this algorithm and repeat the procedure from $t_{k+1}$ to $t_{f}$. 
    \EndIf 
    \State Propagate the RAAN using Eq. \eqref{RAAN1} and \eqref{RAAN2}
    \begin{align} 
  \dot{\Omega}&=-\frac{3}{2} J_2 \sqrt{\frac{\mu}{a^3}} \left(\frac{R_e}{a} \right)^2 \cos{i} \label{RAAN1}\\ 
    \Omega_{t_{k+1}} &=  \Omega_{t_{k}} +  \dot{\Omega}(t_{k+1}-t{k}) \label{RAAN2}
\end{align}
\EndFor
\end{algorithmic}
\end{algorithm}

\subsubsection{RAAN Matching Method}
\begin{figure}[H]
    \centering
    \includegraphics[width = \textwidth]{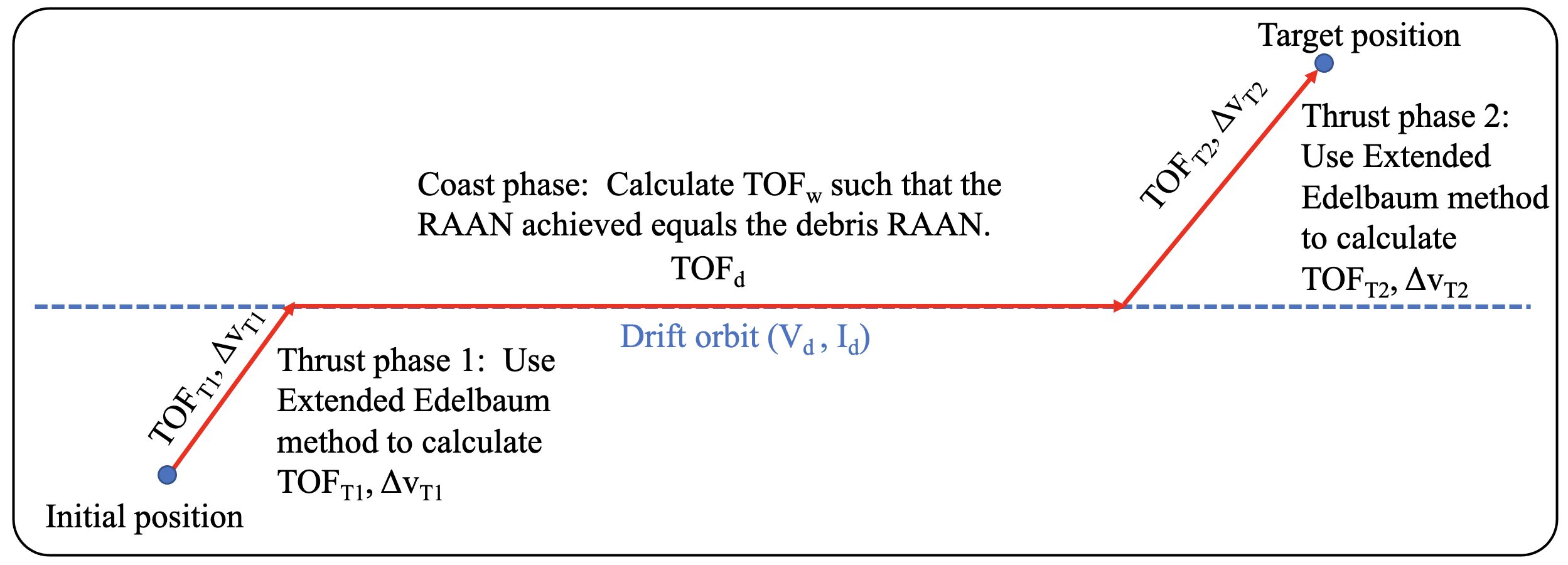}
    \caption{Thrust-drift-thrust structure (subscript $u_1$ shows the first thrust phase, subscript $d$ shows drift phase and subscript $u_2$ indicates the second thrust phase.)}
    \label{f50}
\end{figure}

This method builds on the Extended Edelbaum method such that RAAN changing transfers can be optimized. In this method, orbital precession is used to achieve a target RAAN by drifting at an intermediate drift orbit as done in \cite{6}. The spacecraft follows the thrust-drift-thrust trajectory shown in \Cref{f50} to reach its target state. The drift orbit variables ($V_d$ and $I_d$) are obtained by optimizing the transfer for optimal time or propellant consumption as discussed in Section \ref{optim}.  \par 

During the drifting phase, some thrust may be needed to maintain the orbital altitude. Hence, the thrust magnitude is set to be equal to the drag acceleration experienced, acting in the opposite direction. The two thrust phases in this sequence are evaluated using the Extended Edelbaum method discussed above. \Cref{alg:3} shows the steps associated with calculating the time of flight and $\Delta v$ for a trajectory that utilizes the RAAN matching method to make RAAN changes.

\begin{algorithm}[H]
\caption{RAAN Matching Method}\label{alg:3}
\begin{algorithmic}
\Require Initial orbit elements ($a_0$, $i_0$, $\Omega_{t_0,initial}$), target orbit elements ($a_f$, $i_f$, $\Omega_{t_0,target}$).
\State \textbf{Thrust phase 1}
\State Calculate $TOF_{T1}$ and $\Delta v_{T1}$ using the Extended Edelbaum method. Calculate $\Delta \Omega_{T1}$  (RAAN change of the spacecraft due to precession during thrust phase 1.)
\State  \textbf{Thrust phase 2}
\State Calculate  $TOF_{T2}$ and $\Delta v_{T2}$ using the Extended Edelbaum method. Calculate $\Delta \Omega_{T2}$  (RAAN change of the spacecraft due to precession during thrust phase 2.)
\State \textbf{Drift Phase}
\State Calculate the drift rate of the spacecraft $\Omega_{s/c}$, and the drift rate of the target $\dot{\Omega}_{target}$ using Eq. \eqref{RAAN1}. Then calculate $TOF_d$ (Drift time required to match with the final RAAN), using Eq. \eqref{tw}, which equates the RAAN reached by the Servicer to the RAAN of the debris at arrival time.
\begin{equation}\label{tw}
    (\Omega_{t_0,initial}+ \Delta \Omega_{T1} +\Delta \Omega_{T1} + \dot{\Omega}_{s/c}) TOF_d = \Omega_{t_0, target}+ \dot{\Omega}_{target}(TOF_d+TOF_{T1}+TOF_{T2})
\end{equation}
\State Calculate the $\Delta v$ used to offset drag in the drift phase ($\Delta v_p$). This is achieved by setting the thrust magnitude equal to the drag acceleration (Eq \eqref{eq7}) acting in the opposite direction during drifting. Then,  
\begin{equation}
    \Delta v_p = -\int^{TOF_d}_0 a_{drag} \ dt 
\end{equation}

\State \textbf{Output}: Calculate total $\Delta v$ and $TOF$.
\begin{equation}
    \Delta v=\Delta v_{T1}+ \Delta v_{T2} + \Delta v_p  \ \text{and} \  TOF=TOF_{T1}+TOF_{d}+TOF_{T2}    
\end{equation}
\end{algorithmic}
\end{algorithm}

\subsubsection{Optimization}\label{optim}
The drift orbit parameters and the launch epoch of the mission need to be optimized to obtain the best TOF or $\Delta v$ for a given tour. The input parameters required for the optimization are the duty ratio, maximum thrust, specific impulse, constraints on TOF or $\Delta v$, coordinates of the debris to be removed, Servicer wet mass, and optimization parameter (TOF or $\Delta v$).\par 
The optimization vector ($\bm{x}$) is set to represent all drift orbits in the tour. Thus, for a sequence of N debris, 
\begin{equation}
    \bm{x} = [V_{d_1}, I_{d_1} , V_{d_2} , I_{d_2}, \ldots,
    V_{d_N}, I_{d_N}]
\end{equation}

\par
\begin{figure}[hbt!]
    \centering
    \includegraphics[width = \textwidth]{ 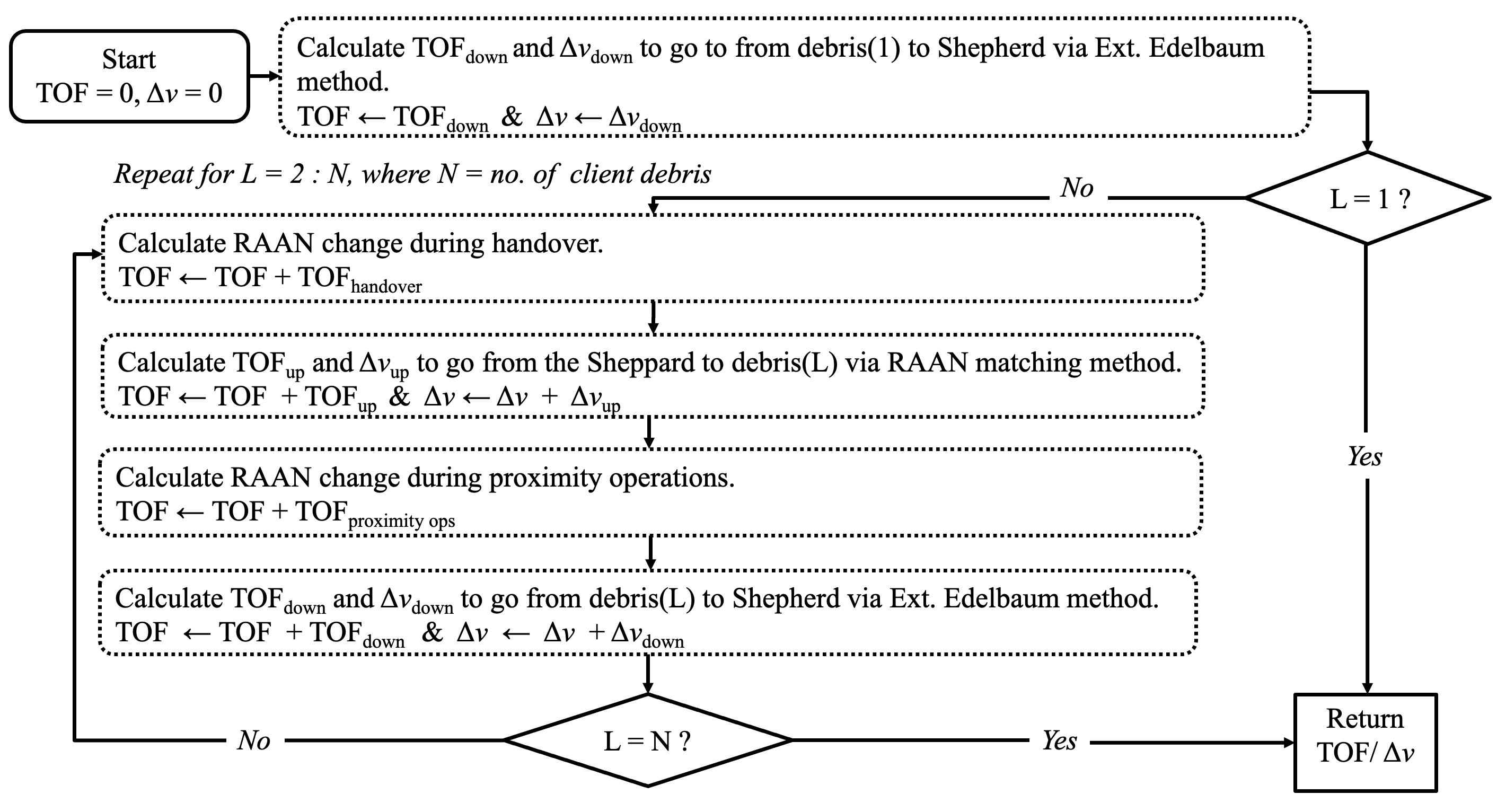}
    \caption{Fitness function calculation for the optimization}
    \label{opt}
\end{figure}

 
Note that the second subscript indicates the debris that is reached by using that drift orbit (i.e. $[V_{d_1}, I_{d_1}]$ is used to reach debris 1 and so on.  )
\Cref{opt} shows the steps involved in the objective function calculation for removing N debris. The optimization uses the interior-point algorithm in Matlab's \textit{fmincon}. As either TOF or $\Delta v$ is being optimized in a given simulation, constraints can be set on the other parameter such that it remains within a feasible domain. These constraints are introduced as nonlinear inequality constraints to the optimization.

\subsection{Mission Guidance}
The PMDT can only provide a low-accuracy solution to the optimal trajectories as the integration of the dynamics is not conducted. If the thrust law computed with the PMDT is applied in a forward manner, the trajectory can deviate significantly from the reference, due to the simplification made in the model and the lack of a feedback mechanism. 
These drawbacks can be circumvented by the inclusion of a guidance scheme that computes the thrust law to track the PMDT reference trajectory. The guidance scheme enables an accuracy assessment of the PMDT models and provides a mean to actually fly the mission. 

Three such schemes are explored to this end: the first one adapts the guidance laws by A. Ruggiero in \cite{e7}, the second
utilizes the $\Delta v$-Law proposed by S. Locoche in \cite{slim}, and the third one 
 uses the Q-law proposed by Petropolous~in~\cite{Petropoulos2005}. Importantly these laws take the transfers computed from PMDT as references to track, thus implicitly exploiting $J_2$: a feature not directly available in these schemes. 
 
\subsubsection{Ruggiero Guidance}

The Ruggiero Guidance was proposed by A. Ruggiero in \cite{e7}. It uses closed loop guidance laws to steer a given orbital element to a target value. To do so, the thrust direction is changed according to the orbital element correction laws generated based on the optimal thrust direction $\boldsymbol{T}$ given in Table 1 in \cite{e7}. Denoting the current orbital elements as $\bm{X} = [a,e,i,\Omega]^\text{T}$ and the target elements as  $\bm{X}_{T} = [a_{T},e_{T},i_{T},\Omega_{T}]^\text{T}$, the optimal thrust vectors for changing each orbital element are: 

\begin{itemize}
    \item For changing semi-major axis (${a}$):
    \begin{equation}\label{a}
    \begin{aligned}
    \bm{T}_a &= (\eta_{a} > \eta_{a}^{bd}) \textrm{sgn}(a_{T} - a)[\cos{\beta_{a}}\sin{\alpha_{a}}, \cos{\beta_{a}}\cos{\alpha_{a}}, \sin{\beta_{a}}] \\ 
    &\text{where}
    \tan{\alpha_{a}} = \frac{{e} \sin{{\nu}}}{1+ {e} \cos{{\nu}}}, \beta_{a} = 0, \eta_{a} = |\bm{v}| \sqrt{\frac{a (1- {e})}{\mu(1+ {e})}}
        \end{aligned}
    \end{equation}
    \item For changing eccentricity ($e$):
      \begin{equation}\label{e}
    \begin{aligned}
    \bm{T}_e &= (\eta_e > \eta_e^{bd}) sgn(e_{T} - {e})[\cos{\beta_{e}}\sin{\alpha_{e}}, \cos{\beta_{e}}\cos{\alpha_{e}}, \sin{\beta_{e}}] \\ 
    &\text{where}
    \tan{\alpha_e} = \frac{\sin{{\nu}}}{\cos{{E}} +\cos{{\nu}}}, \beta_e = 0, \eta_e = \frac{1+ 2{e} \cos{{\nu} + \cos^2{{\nu}}}}{1+ {e} \cos{{\nu}}}
        \end{aligned}
    \end{equation}
    \item For changing inclination ($i$):
          \begin{equation}\label{i}
    \begin{aligned}
    \bm{T}_i &= (\eta_i > \eta_i^{bd}) \textrm{sgn}(i_{T} - {i})[\cos{\beta_{i}}\sin{\alpha_{i}}, \cos{\beta_{i}}\cos{\alpha_{i}}, \sin{\beta_{i}}] \\ 
    &\text{where}
    \tan{\alpha_i} = 0,  \beta_i = \frac{\pi}{2} \textrm{sgn}(\cos({\omega +\nu})), \eta_i  =\frac{|\cos (v+\omega)|}{1+e \cos v}\left(\sqrt{1-\mathrm{e}^{2} \sin ^{2} \omega}-\mathrm{e}|\omega|\right)
        \end{aligned}
    \end{equation}
\item For changing the right ascension of the ascending node ($\Omega$):
             \begin{equation}\label{omega}
    \begin{aligned}
    \bm{T}_\Omega &=  (\eta_\Omega > \eta_\Omega^{bd}) \textrm{sgn} \left(\frac{-\sin{(\Omega - \Omega_{T})}}{\sqrt{1 - \cos^2{(\Omega -\Omega_{T})}}}\right)[\cos{\beta_{\Omega}}\sin{\alpha_{\Omega}}, \cos{\beta_{\Omega}}\cos{\alpha_{\Omega}}, \sin{\beta_{\Omega}}] \\ 
    &\text{where}
    \tan{\alpha_\Omega} = 0,  \beta_\Omega = \frac{\pi}{2} \textrm{sgn}(\sin({\omega +\nu})), \eta_\Omega = \frac{|\sin{(\nu+\omega)}|}{1+e\cos {\nu}}(\sqrt{1- e^2\cos^2\nu } - e|\sin\omega|)
        \end{aligned}
    \end{equation}
(Note that the term $\textrm{sgn} \left(\frac{-\sin{(\Omega - \Omega_{T})}}{\sqrt{1 - \cos^2{(\Omega -\Omega_{T})}}}\right)$ is used instead of $\textrm{sgn}(\Omega_T - \Omega)$ to find the direction of the thrust required to go towards $\Omega_T$. This is done to find the correct direction with the smallest angle between $\Omega$ and $\Omega_T$)
\end{itemize}

 $\eta^{bd}$ is a limit beyond which the thrust vector shall be activated. At each time step of the tour, the orbital elements from the PMDT solution are taken as the target when providing guidance.   \par
Once the optimal thrust vectors are calculated, weighting coefficients ($c_X$) are introduced in front of each. Then, the optimal unit thrust acceleration ($\bm{u}$)  can be calculated as follows. 
\begin{equation} \label{19}
    \bm{u} = \frac{c_a \bm{T}_a + c_e \bm{T}_e + c_i \bm{T}_i + c_{\Omega} \bm{T}_{\Omega}}{|c_a \bm{T}_a + c_e \bm{T}_e + c_i \bm{T}_i + c_{\Omega} \bm{T}_{\Omega}|} \ \text{where} \ c_X = {|X- X_{T}|} W_x
\end{equation}
For each element in $\bm{X}$, $c_X$ is dependent on the difference between the current value of the orbital element and the target value. $\bm{W} =  [W_a,W_e,W_i,W_\Omega]^\text{T}$ are coefficients to be optimized to improve the performance of Ruggiero guidance.

\subsubsection{{$\Delta v$}-Law}

The guidance scheme called the $\Delta v$-Law proposed by Slim Locoche in \cite{slim} is based on Lyapunov feedback control. 
This method entails developing a control feedback algorithm that decreases a scalar function (called the Lyapunov function), representing the distance between the current state and its target. The designed control algorithm aims to drive the  Lyapunov function ($L$) to zero. The $L$ function used in the $\Delta v$-Law is
\begin{equation}
L=\widetilde{\Delta v}^2 \equiv \lambda_a\left[V_c^2-2 V_c V_{c f} \cos (\pi / 2 \Delta \sigma)+V_{c f}^2\right]+\frac{4}{9} \lambda_{e_1}\left[\frac{\left[\left(1-\lambda_{e_2}\right) V_c+\lambda_{e_2} V_{c f}\right]\left[\operatorname{asin}(e)-\operatorname{asin}\left(e_f\right)\right]}{\cos (\widetilde{\beta})}\right]^2
\end{equation}

where:
\begin{equation}
\begin{gathered}
\Delta \sigma=\sqrt{\left[\lambda_{a, i} \Delta i\right]^2+\left[\lambda_{a, \Omega} \sin (\mathrm{i}) \Delta \Omega\right]^2} \\
\tan (|\tilde{\beta}|)=\left|\frac{3 \pi \lambda_{e, i} \Delta i}{4 \cos \left(\lambda_\omega \omega\right)\left[\ln \left(\frac{e_f+1}{e_f-1}\right)+\ln \left(\frac{e-1}{e+1}\right)-\Delta e\right]}\right|
\end{gathered}
\end{equation}

Here, $V_c$ and $V_f$ indicate the current and target orbital velocities. $\lambda_{e_1}, \lambda_{e_2}, \lambda_{a, i}, \lambda_{e, i}, \lambda_{a, \Omega}, \lambda_\omega$ are parameters to be optimized to enhance the performance of the Lyapunov controller developed. Note that $L$ is constructed by combining the analytical $\Delta v$ equations for making semi-major axis, eccentricity, inclination, and RAAN changes \cite{slim}. Also note that $\Delta i = i - i_T$, $\Delta \Omega = \cos^{-1}( \cos(\Omega - \Omega_T))$ and $\Delta e = e - e_T$. \par 
The goal of Lyapunov control is to make $ \dot{L}$ as negative as possible, such that $L$ shall approach zero quickly. Note that \begin{equation}\label{equ:Ldot}
\begin{aligned}
    \dot{L} = \frac{\partial L}{\partial \bm{X}} \dot{\bm{X}} = \frac{\partial L}{\partial \bm{X}} \bm{B} \bm{u},
\end{aligned}
\end{equation}
where $\bm{X}$ denotes the state variables as given in Ruggiero guidance and 
$\bm{B}$ represents the Gauss Variational Equations (GVEs) for the slow variables given by
\def\GVEexp{\begin{matrix}
    \frac{2a^2}{h} e\sin(\nu) & \frac{2a^2}{h} \frac{p}{r} & 0 \\
    \frac{1}{h} p\sin(\nu) & \frac{1}{h} ((p+r)\cos(\nu)+re)  & 0 \\
    0 & 0 & \frac{r\cos(\omega + \nu)}{h} \\
    0 & 0 &\frac{r\sin(\omega + \nu)}{h\sin(i)} \\
\end{matrix}}
\begin{equation}\label{equ:LagrangeEquations}
    \bm{B} = \left[ \GVEexp \right].
\end{equation}
Hence, the optimal control acceleration direction can be calculated as 
\begin{equation}\label{ulopt}
\begin{aligned}
    \bm{u} = - \frac{\bm{B}^\text{T} \left( \frac{\partial L}{\partial \bm{X}} \right)^\text{T}}{\left|\left| \left( \frac{\partial L}{\partial \bm{X}} \right) \bm{B} \right|\right|}
\end{aligned}
\end{equation}




\subsubsection{Q-Law}
One of the most versatile and well-known control laws is the Q-Law developed by Petropoulos \cite{Petropoulos2005}. The Q-law is best thought of as a weighted, squared summation of the time required to change the current state $\bm{X} = [a,e,i,\Omega]^\text{T}$ to the target state $\bm{X}_T=[a_T,e_T,i_T,\Omega_T]^\text{T}$. It can be written as 
\begin{equation}\label{equ:Qlaw}
\begin{aligned}
    Q = (1+W_P P(\bm{X}))\sum_{X} S_X(\bm{X}) W_X(\bm{X}) \left( \frac{\delta (X,X_T)}{\text{max}_\nu (\dot{X})} \right)^2,
\end{aligned}
\end{equation}
where $W_P$ and $P$ form a penalty function and $S_{X}$ are scaling functions. These are functions of the state and can be found in \cite{Petropoulos2005}. $\delta(X,X_T) = {X} - {X}_{T}$ for ${X} = a,e,i$ whilst $\delta(X,X_T) = \cos^{-1}(\cos({X} - {X}_{T}))$ for ${X} = \Omega$. The expressions $\text{max}_\nu (\dot{{X}})$ are the maximum rate of change of each COE over the current osculating orbit and can be calculated analytically for all elements. The weights $W_X$ can be used to prioritize which elements to target.

Via Lyapunov's second theorem, a stable control is one that ensures $\dot{Q}<0$ throughout the transfer. One way of doing this is to select a controller that minimizes the rate of change of the Lyapunov function (in this case, the most negative value).
\begin{equation}\label{equ:Qdot}
\begin{aligned}
    \dot{Q} = \frac{\partial Q}{\partial \bm{X}} \dot{\bm{X}} = \frac{\partial Q}{\partial \bm{X}} \bm{B} \bm{u},
\end{aligned}
\end{equation}

Leading to a control acceleration direction
\begin{equation}\label{equ:Uoriginal}
\begin{aligned}
    \bm{u} = - \frac{\bm{B}^\text{T} \left( \frac{\partial Q}{\partial \bm{X}} \right)^\text{T}}{\left|\left| \left( \frac{\partial Q}{\partial \bm{X}} \right) \bm{B} \right|\right|}
\end{aligned}
\end{equation}

\subsubsection{Propagating with Guidance}

\paragraph{Integrating the effect of eclipses and duty ratio} \par 

\begin{figure}[hbt!]
    \centering
    \includegraphics[width = 0.6\textwidth]{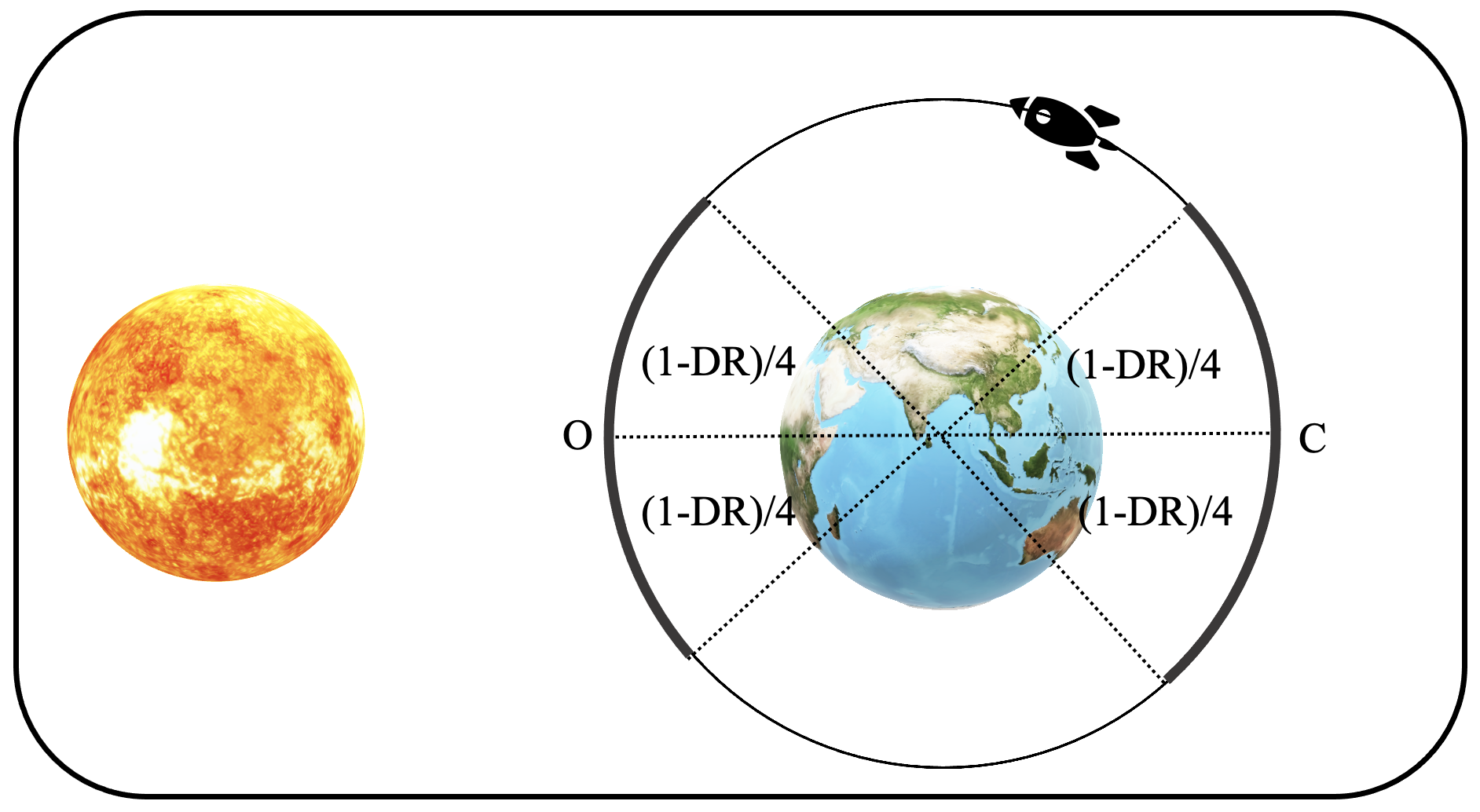}
    \caption{Eclipse formulation for propagating with guidance}
    \label{eclprop}
\end{figure}

When propagating the dynamics with guidance, the effect of eclipses and the duty ratio must be considered adequately. 
 However, turning thrust off asymmetrically (i.e., only during the eclipse) will cause eccentricity buildup \cite{Viavattene2022DesignOM}. Hence thrust is turned off symmetrically across the orbit in the highlighted regions in Figure \ref{eclprop}. These regions are determined as follows.\par 
First, the value of the argument of latitude at the eclipse center ($\theta_C$) is calculated at each time step of the propagation. Then, the thrust is turned off 
for a symmetric fraction of the orbit when 
$  \theta_C - 2 \pi  \frac{1-DR}{4}\leq \theta \leq \theta_C + 2 \pi\frac{1- DR}{4}$, where $DR$ is the duty ratio. Then, the thrust is also turned off for the same orbit fraction opposite from C (at O), when $\theta_O- 2 \pi\frac{1- DR}{4}\leq \theta \leq \theta_O +2 \pi \frac{1- DR}{4}$

\paragraph{Counteracting drag in drift orbits}\par  \label{drag}
It was noted that when using guidance in the drift orbits used for RAAN matching, a large amount of propellant shall be used to counteract minute orbital changes made by the drag. To minimize this, the thrust was only turned on when the drag resulted in significant orbit changes compared to the reference drift orbit. Hence, on drift orbits, the thrust was only turned on-provided that the spacecraft was not in eclipse- when $|a - a_{ref}| > 5$ km or $|i - i_{ref}| > 0.1$ deg or $|\Omega - \Omega_{ref}| > 0.1$ deg. The thrust was turned off when $|a - a_{ref}| < 0.5$ km and $|i - i_{ref}| < 0.01$ deg and $|\Omega - \Omega_{ref}| < 0.01$ deg, to conserve fuel. \par 

Algorithm \ref{alg:5} illustrates how the three guidance schemes are implemented when propagating the dynamics. Note that the mean orbital elements are used to calculate the control acceleration direction through guidance, such that propellant is not wasted on correcting for the osculating nature of the elements.

\begin{algorithm}[hbt!]
\caption{Implementation of Guidance}\label{alg:5}
\begin{algorithmic}
\setstretch{1}
\Require The PMDT solution of the transfer, the transfer time ($t_f$). 
\For {$t = 0 : t_f$} 
\State 1. Convert the osculating orbital elements to mean elements.
\State 2. Calculate the target state ($\bm{X}_T$) at $t$ by interpolating the PMDT solution. 
\State 3. Calculate the optimal thrust direction ($\bm{u}$) using Eq \eqref{19} (For Ruggiero guidance), Eq. \eqref{ulopt} (for $\Delta v$-Law guidance) or Eq. \eqref{equ:Uoriginal} (for Q-Law guidance), taking the mean elements as the current state.
\State 4. Calculate the effect of the eclipse and duty ratio.
\setlength\parindent{24pt}
\If{ ($\theta_C - 2 \pi  \frac{1-DR}{4}\leq \theta \leq \theta_C + 2 \pi\frac{1-DR}{4}$ \textbf{or}  $\theta_O- 2 \pi\frac{1-DR}{4}\leq \theta \leq \theta_O +2 \pi \frac{1-DR}{4}$)} 
$\eta=0$
\Else{}
$\eta=1$
\EndIf
\State 5. If propagating a drift orbit, implement the drag counteraction method.
\If{( $|a - a_{ref}| > 5$ km \textbf{or} $|i - i_{ref}| > 0.1$ deg \textbf{or} $|\Omega - \Omega_{ref}| > 0.1$ \textbf{and} $\eta \neq 0$)} 
$\eta=1$
\ElsIf{($|a - a_{ref}| < 0.5$ km \textbf{and} $|i - i_{ref}| < 0.01$ deg \textbf{and} $|\Omega - \Omega_{ref}| < 0.01$ deg)}
$\eta=0$
\EndIf
\State 6. Calculate thrust acceleration ($\bm{a}_{T}$), where $m$ is the system mass.
\begin{equation}
   \bm{a}_{T}   = \eta \frac{f}{m}\bm{u}
\end{equation}
\State 7. Calculate the $J_2$ acceleration ($\bm{a}_{J_2}$) using the formulation shown in \cite{8}.

\State 8. Calculate the drag acceleration ($\bm{a}_{drag}$), as described in \Cref{alg:1}. 
\State 9. Calculate the time derivatives of state and spacecraft mass. 
\begin{align}
        \ddot{\bm{X}} &= \bm{a}_{J2} + \bm{a}_{drag} + \bm{a}_{T}+\bm{a}_{gravity}\\
    \dot{m} &= \frac{T_{max}}{I_{sp} g_0}
\end{align}

\EndFor 
\State Propagate the dynamics to obtain $\bm{X}$ using a propagator such as Matlab's \textit{ode45}. 
\end{algorithmic}
\end{algorithm}

\subsubsection{Optimizing the guidance parameters}
Each of the three guidance approaches has user-defined parameters in the form of weights: $W_x$ for Ruggiero and the Q-Law guidance, and $\lambda_x$ for the $\Delta v$-Law guidance. These parameters can significantly affect the behavior of the guidance control laws. As such, a Particle Swarm Optimisation (PSO) is used to tune these parameters. Due to the computational challenges of converting from osculating to mean elements and computing the duty cycle, the guidance strategy was simplified for the PSO.  \par 
A duty cycle of $50\%$ was simulated by reducing the spacecraft thrust by a factor of~$2$. In addition, only the average/secular contribution of $J_2$ was considered, removing the need for converting from osculating to mean elements when computing the control. Given that there are two distinct transfer scenarios, either lowering the altitude with the debris attached to the Servicer or climbing and rendezvousing to the debris, different sets of weights were computed for these two transfer scenarios. This further reduced the computational requirements and created a design vector of 10 elements, 5 for the downward legs and 5 for the upwards legs. In the case of the Ruggiero and Q-law guidance, the weights relating to the eccentricity and argument of periapsis can automatically be set to $0$. In the case of Locoche guidance, the problem can be simplified by assuming $\lambda_\omega=0$ as the eccentricity and the argument of periapsis are not targeted. 
\par 
Six PSO simulations were initiated, with a swarm size of 50 each, tracking the time- and fuel-optimal reference trajectories with the Ruggiero, $\Delta v$-Law, and Q-Law guidance strategies. Whilst tracking the reference trajectory, the objective is to minimise the accumulated errors in semi-major axis, inclination and RAAN at the end of each leg. Once a set of weights were obtained in the simplified scenarios used in the PSO, there were deployed in the original, osculating dynamics with the $50\%$ duty cycle. The coefficients provided by these PSO simulations are reported in Table~\ref{table:PSOcoeff}. These coefficients were subsequently used for the guidance reported in Sections~\ref{sec:FuelResults} and~\ref{sec:TimeResults}. Naturally, there is a deviation from the PSO results, but the guidance remains satisfactory without requiring extensive computational resources. 

\begin{table}[hbt!]
\centering 
\caption{Guidance coefficients obtained from PSO simulations in simplified dynamics}
\label{table:PSOcoeff}
\begin{tabular}{c|ccccc|ccccc}
\toprule
Guidance       & \multicolumn{5}{c|}{Time Optimal}                                                       & \multicolumn{5}{c}{Fuel Optimal}                                                        \\ \hline
Ruggiero       & $W_a$          & $W_e$          & $W_i$          & $W_\Omega$     & -                   & $W_a$          & $W_e$          & $W_i$          & $W_\Omega$     & -                   \\ \hline
Leg 1          & 0.2058         & 0              & 0              & 0              & -                   & 0.2058         & 0              & 0              & 0              & -                   \\
Leg 2          & 0.9564         & 0              & 0.1211         & 0.00985        & -                   & 0.8622         & 0              & 1              & 0              & -                   \\
Leg 3          & 0.2058         & 0              & 0              & 0              & -                   & 0.2058         & 0              & 0              & 0              & -                   \\
Leg 4          & 0.9564         & 0              & 0.1211         & 0.00985        & -                   & 0.8622         & 0              & 1              & 0              & -                   \\
Leg 5          & 0.2058         & 0              & 0              & 0              & -                   & 0.2058         & 0              & 0              & 0              & -                   \\ \hline
$\Delta v$-Law & $\lambda_{e1}$ & $\lambda_{e2}$ & $\lambda_{ai}$ & $\lambda_{ei}$ & $\lambda_{a\Omega}$ & $\lambda_{e1}$ & $\lambda_{e2}$ & $\lambda_{ai}$ & $\lambda_{ei}$ & $\lambda_{a\Omega}$ \\ \hline
Leg 1          & 0.0169         & 0.9431         & 0.9998         & 0.9010         & 0.00329             & 0.0341         & 0.9595         & 0.7910         & 0              & 0.0006800           \\
Leg 2          & 0.2312         & 0.4254         & 1              & 0.1943         & 0.0724              & 0.1260         & 0.04625        & 0.5969         & 0.8367         & 0.01846             \\
Leg 3          & 0.0169         & 0.9431         & 0.9998         & 0.9010         & 0.00329             & 0.0341         & 0.9595         & 0.7910         & 0              & 0.0006800           \\
Leg 4          & 0.2312         & 0.4254         & 1              & 0.1943         & 0.0724              & 0.1260         & 0.04625        & 0.5969         & 0.8367         & 0.01846             \\
Leg 5          & 0.0169         & 0.9431         & 0.9998         & 0.9010         & 0.00329             & 0.0341         & 0.9595         & 0.7910         & 0              & 0.0006800           \\ \hline
Q-Law          & $W_a$          & $W_e$          & $W_i$          & $W_\Omega$     & -                   & $W_a$          & $W_e$          & $W_i$          & $W_\Omega$     & -                   \\ \hline
Leg 1          & 0.8784         & 0              & 0.9156         & 0.0001773      & -                   & 1.0              & 0              & 0.8479         & 0              & -                   \\
Leg 2          & 1.0            & 0              & 0.04448        & 0.009775       & -                   & 0.6934         & 0              & 0.5649         & 0              & -                   \\
Leg 3          & 0.8784         & 0              & 0.9156         & 0.0001773      & -                   & 1.0               & 0              & 0.8479         & 0              & -                   \\
Leg 4          & 1.0            & 0              & 0.04448        & 0.009775       & -                   & 0.6934         & 0              & 0.5649         & 0              & -                   \\
Leg 5          & 0.8784         & 0              & 0.9156         & 0.0001773      & -                   & 1.0               & 0              & 0.8479         & 0              & -                   \\ \bottomrule
\end{tabular}
\end{table}

\section{Exemplar 3-Debris ADR Mission} \label{results}
This section discusses the results of optimizing an ADR mission for three debris in near sun-synchronous orbits. The objects to be removed are, in order, H-2A R/B (ID: 33500), ALOS 2 (ID: 39766), and GOSAT (ID: 33492). This sequence of objects has been chosen arbitrarily to provide an exemplar test case on which PMDT could be evaluated. Note that the debris sequence could be optimized using the PMDT to asses various sequences, but this was not performed in this preliminary study.  \par 
The propulsion conditions used were 60~mN maximum thrust, 50\% duty ratio, and 1300 s specific impulse. The Servicer was assumed to have a wet mass of 800~kg, and the launch date was set to be 25-Mar-2022 06:37:09 UTC. The altitude at which the debris and the Servicer meet the Shepherd was set to 350~km, i.e., below the ISS altitude, to satisfy safety requirements. Eclipses and drag acceleration were taken into consideration for this case study. Once the Servicer has reached the debris,  45~days were allocated to perform proximity operations before the start of deorbiting. When the Servicer and the debris reached the 350~km circular orbit, 30~days were allocated for the handover of the debris from the Servicer to the Shepherd. 

\subsection{\label{sec:FuelResults}Fuel-Optimal Scenario}
When generating the fuel-optimal result, the total TOF was limited to less than five years (1825~days).
In the optimized solution, the first drift orbit (going from the Shepherd altitude to ALOS 2) was at $a = 7662.8$~km and $ i = 98.29$ deg. The second optimal drift orbit (going from the Shepherd altitude to GOSAT) was at $a = 7499.2$~km and $ i =  97.85$ deg. The optimal $\Delta v $ obtained was 945.58 m/s, with a total fuel consumption of 136.07~kg. Note that the optimized tour requires five years, i.e., it is at the limit for the total allowed TOF. Tours with better fuel consumption can be obtained at the cost of increased mission duration. 
\par
\begin{figure}[hbt!]
    \centering
    \includegraphics[width = 0.8 \textwidth]{ 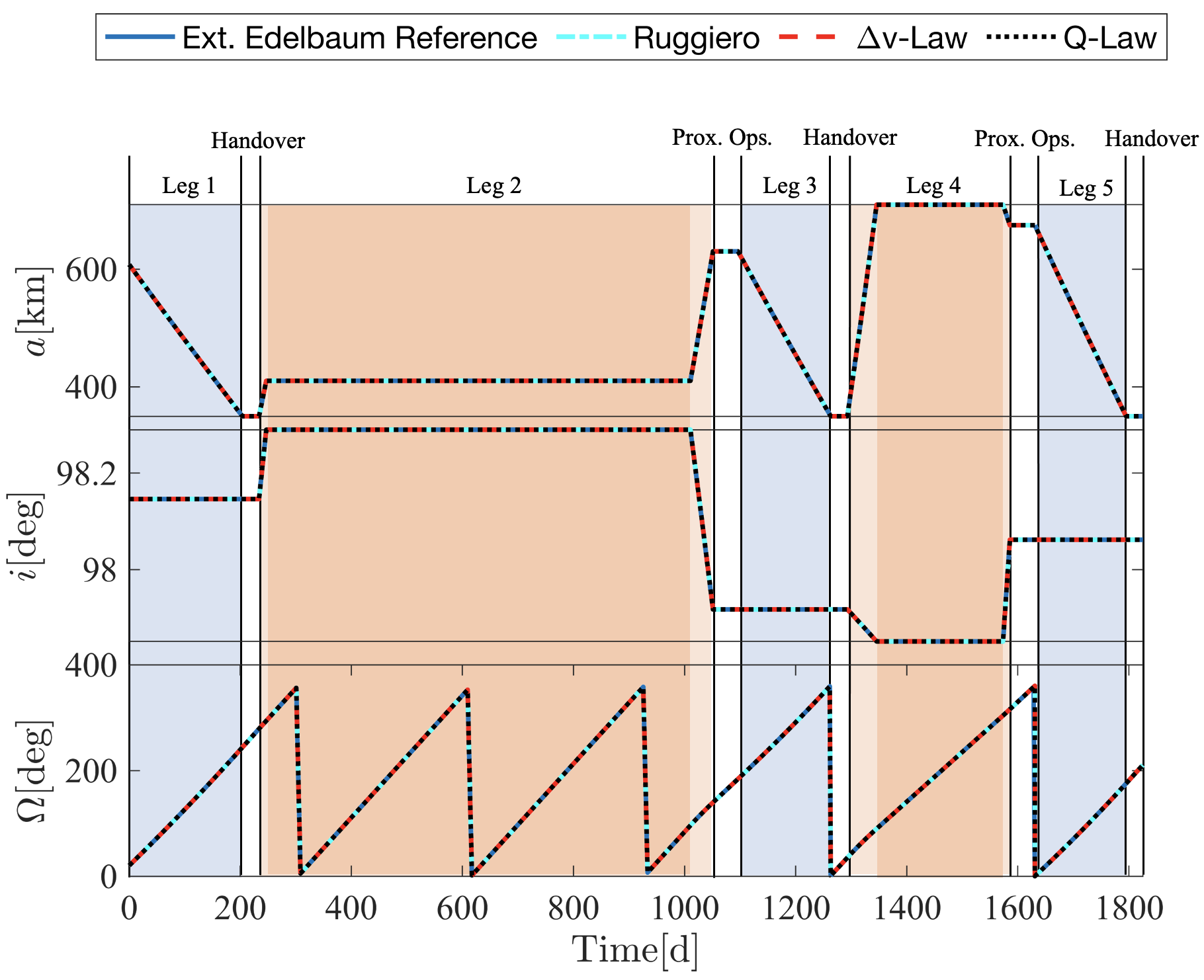}
    \caption{Fuel optimal tour with guidance}
    \label{l123}
\end{figure}

 
\Cref{t8} shows the $\Delta$v and TOF breakdown of the fuel optimal trajectory obtained. \Cref{l123}  show the plots of semi-major axis, inclination, and RAAN variations observed during the tour. \par 
\Cref{t3,error_fopt_table} show the performance comparison between the three guidance schemes and the forward propagation of the PMDT outcome.  The forward propagation of the drift orbits are done by assuming that at each time step, $\bm{a}_T = -\eta \frac{\bm{a}_{Drag}}{DR}$. For all other legs, the out of plane thrust angle ($\beta$) obtained from the Extended Edelbaum method is interpolated at each time step to calculate $\bm{a}_T = \eta \frac{f}{m}[0, \cos{(\beta_{interp})}, \sin{(\beta_{interp})}]$.\par 
 Note that it is assumed that errors that are smaller than $\approx$20 km in semi-major axis, $\approx$1 deg in RAAN, and $\approx$ 0.1 deg in inclination shall be taken care of in the proximity operation and handover phases of the mission. The coefficients of guidance were optimized by PSO such that the errors at each leg would not exceed these limits. It can be seen that the propagation errors are reduced through the use of either guidance scheme.  \par
\begin{figure}[hbt!]
    \centering
    \includegraphics[width = 0.8 \textwidth]{ 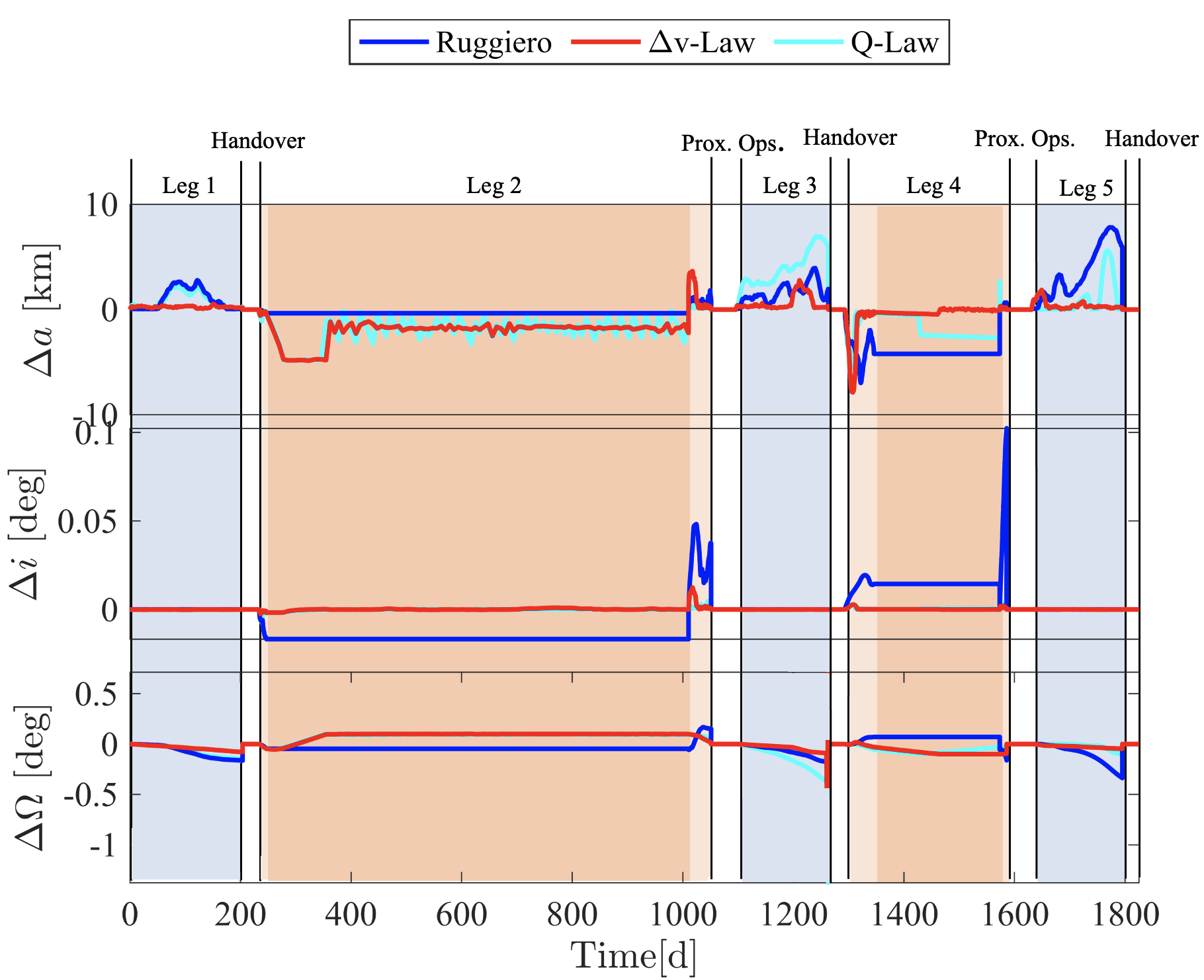}
    \caption{Fuel optimal tour guidance errors present}
    \label{error_fopt}
\end{figure}

It is also evident that the Ruggiero guidance consumes the least fuel. $\Delta v$-Law consumes more propellant but shows better accuracy in tracking semi-major axis and inclination.  The Q-Law consumes even more fuel, but without achieving a better accuracy. As figure \ref{error_fopt} shows, Ruggiero guidance can track the reference RAAN through the drift orbit in Leg 2 to a better extent than the $\Delta v$-Law and Q-Law.


\begin{table}[hbt!]
\centering 
\caption{Fuel optimal tour: $\Delta$v and TOF per leg}
\label{t8}
\begin{tabular}{lll}
\toprule
Leg                                   & $\Delta$v (m/s) & TOF (d) \\ \hline
Leg 1 (from H-2A R/B to 350 km orbit) & 140.12             & 203.82     \\
Handover                              & -               & 30.00      \\
Leg 2 (from 350 km orbit to ALOS 2)   & 238.54             & 817.44     \\
Proximity Operations                  & -               & 45.00      \\
Leg 3 (from ALOS 2 to 350 km orbit)   & 151.93             &166.91    \\
Handover                              & -               & 30.00      \\
Leg 4 (from 350 km orbit to GOSAT)    & 239.86            &292.57    \\
Proximity Operations                  & -               & 45.00      \\
Leg 5 (from GOSAT to 350 km orbit)    & 175.14            & 164.25    \\
Handover                              & -               & 30.00      \\ \hline
Total                                 & 945.58            & 1,825.00 \\   \bottomrule
\end{tabular}
\end{table}

\begin{table}[hbt!]
\centering 
\caption{Fuel optimal tour: Error per leg obtained for each of the guidance schemes used}
\label{error_fopt_table}
\begin{tabular}{c|ccc|ccc} \hline
Leg no. & \multicolumn{3}{c|}{Forward Propagated PMDT}               & \multicolumn{3}{c}{Ruggeiro guidance}                      \\ 
        & $\Delta a$ (km) & $\Delta i$ (deg) & $\Delta \Omega$ (deg) & $\Delta a$ (km) & $\Delta i$ (deg) & $\Delta \Omega$ (deg) \\ \hline
1       & 2.9424          & 0.000044         & 0.3086                & 0.0847          & 0.000004         & 0.1641                \\
2       & 11.5665         & 0.050289         & 1.6496                & 0.8959          & 0.037086         & 0.1514                \\
3       & 7.7504          & 0.000006         & 0.7289                & 1.9466          & 0.000004         & 0.1794                \\
4       & 1.6705          & 0.023349         & 0.2748                & 0.5171          & 0.102374         & 0.1619                \\
5       & 5.7491          & 0.000005         & 0.0060                & 5.9474          & 0.000013         & 0.3375                \\ \hline
Total   & 29.6790         & 0.074            & 2.9680                & 9.3916          & 0.139            & 0.9943                \\ \hline \bottomrule
Leg no. & \multicolumn{3}{c|}{$\Delta v$-Law Guidance}               & \multicolumn{3}{c}{Q-Law Guidance}                         \\ 
        & $\Delta a$ (km) & $\Delta i$ (deg) & $\Delta \Omega$ (deg) & $\Delta a$ (km) & $\Delta i$ (deg) & $\Delta \Omega$ (deg) \\ \hline 
1       & 0.1589          & 0.000068         & 0.0800                & 0.2785          & 0.000003         & 0.1470                \\
2       & 0.3412          & 0.001            & 0.7606                & 0.2185          & 0.000745         & 0.0189                \\
3       & 0.2365          & 0.000062         & 0.0906                & 5.8769          & 0.000003         & 0.3859                \\
4       & 0.1256          & 0.000888         & 0.0957                & 0.0237          & 0.000631         & 0.6547                \\
5       & 0.2256          & 0.00008          & 0.0467                & 0.0125          & 0.000023         & 0.0933                \\ \hline
Total   & 1.0877          & 0.002            & 1.0736                & 6.4102          & 0.001            & 1.2998                \\ \hline \bottomrule 
\end{tabular}
\end{table}

\begin{table}[hbt!]
\centering 
\caption{Fuel optimal tour: $\Delta$v, TOF and total error comparison of guidance strategies}
\label{t3}
\begin{tabular}{cccccc}\toprule
Case                    & TOF   (d) & Fuel   (kg) & $\Delta a$ (km) & $\Delta i$ (deg) & $\Delta \Omega$ (deg) \\ \hline
Forward Propagated PMDT  & 1,825.0    & 136.10      & 29.679 & 	0.074& 	2.968                   \\
Ruggiero guidance       & 1,825.4    & 137.11      & 9.392           & 0.139            & 0.994                 \\
$\Delta v$-Law guidance & 1,824.6    & 138.32      & 1.088           & 0.002            & 1.074                 \\
Q-Law guidance          & 1,825.0    & 138.98      & 6.410           & 0.001            & 1.300                 \\
\bottomrule
\end{tabular}
\end{table}

\subsection{\label{sec:TimeResults}Time Optimal Scenario}
When generating the time-optimal result, the total $\Delta v$ was limited to less than 1,500~m/s. The 
first optimal drift orbit was at $a = 7714.2$ km and $ i = 99.35$ deg. The second optimal drift orbit was at $a = 7,464.5$ km and $ i =  97.60$ deg. The optimal time of flight obtained was 1,274.54~days, which occurred when $\Delta v = 1,500.00 $~m/s and the total fuel consumption was 166.18~kg. As this solution reached the $\Delta v$  boundary, it is evident that outcomes with better TOF can be obtained at the cost of increasing the $\Delta v$.

\par
\begin{figure}[hbt!]
    \centering
    \includegraphics[width = 0.8 \textwidth]{ 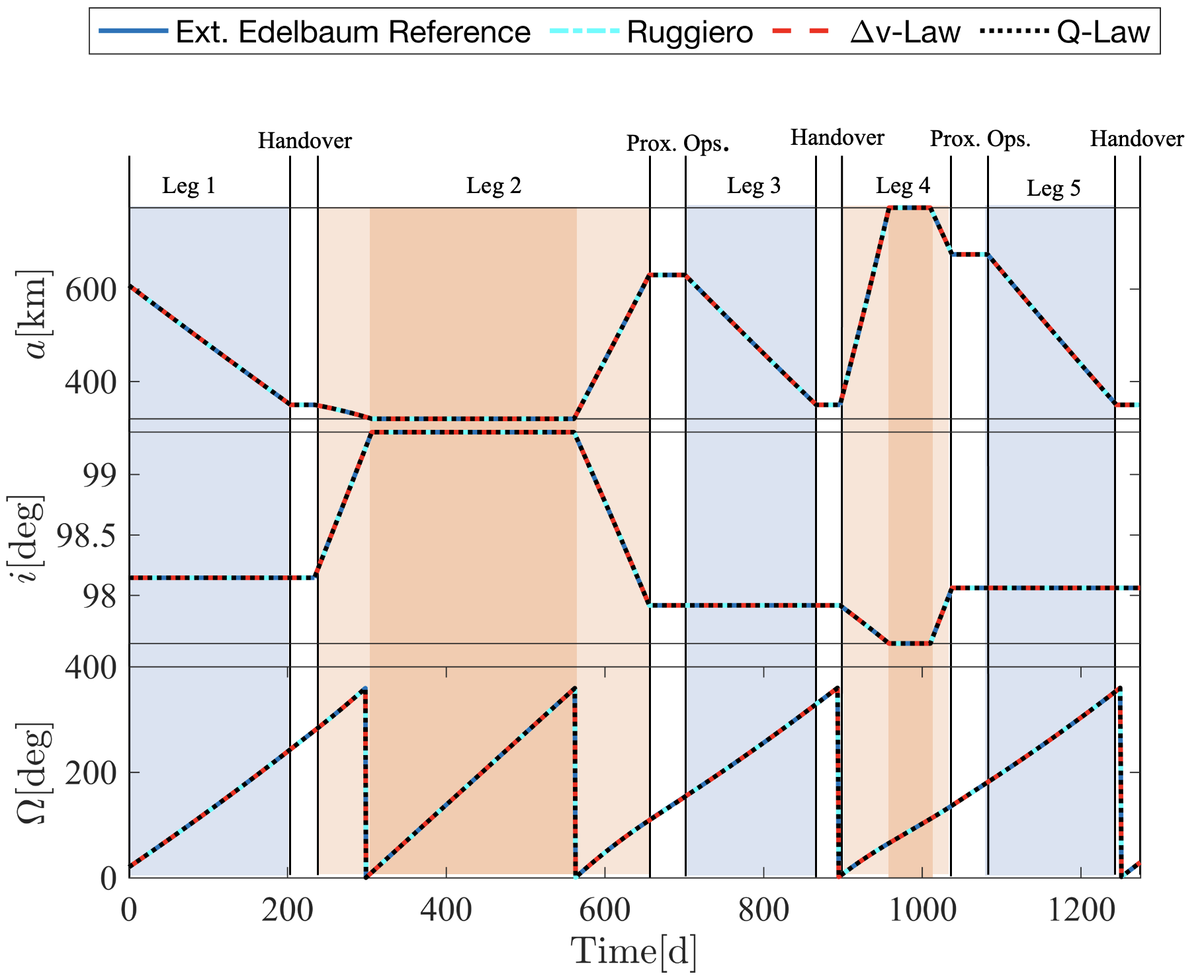}
    \caption{Time optimal tour with guidance}
    \label{t123}
\end{figure}

 
\Cref{t2} shows the $\Delta$v and TOF breakdown of the obtained time-optimal trajectory. \Cref{t123} shows the plots of semi-major axis, inclination, and RAAN variations observed throughout the tour. 
\Cref{t32,error_topt_table}  illustrate how the guidance solution varies from the forward propagated PMDT solution at various mission stages for the time optimal case.  The forward propagation has a significant RAAN and inclination error in the second leg, largely due to the low drift orbit altitude. 
It can be seen that the Ruggiero guidance performs better than the $\Delta v$-Law guidance, but consumes the highest amount of propellant. 
The Q-Law guidance consumes the least amount of propellant and reduces the RAAN error. The $\Delta v$-Law guidance provides a middle ground, where the inclination error is reduced at the cost of increasing the RAAN and semi-major axis errors. 

The results of the time-optimal scenario reinforce the conclusion about the suitability of the PMDT for ADR mission design and of the guidance schemes in tracking the reference solutions.

\begin{table}[hbt!]
\centering 
\caption{ Time optimal tour: $\Delta v$ and TOF per leg}
\label{t2}
\begin{tabular}{lll}
\toprule
Leg       & $\Delta v$ (m/s) & TOF (d) \\\hline 
Leg 1 (from H-2A R/B to 350 km orbit)    & 140.12      & 203.82      \\
Handover   & -        & 30.00        \\
Leg 2 (from 350 km orbit to ALOS 2)     &691.14     & 422.09      \\
Proximity Operations & -        & 45.00        \\
Leg 3 (from ALOS 2 to 350 km orbit)      & 151.93      &165.40      \\
Handover  & -        & 30.00        \\
Leg 4 (from 350 km orbit to GOSAT) & 341.67     & 141.03      \\
Proximity Operations  & -        & 45.00        \\
Leg 5 (from GOSAT to 350 km orbit)      & 175.14      & 162.20     \\
Handover  & -        & 30.00       \\ \hline
Total  & 1,500.00 & 1,274.54 \\ 
\bottomrule
\end{tabular}
\end{table}

\begin{figure}[hbt!]
    \centering
    \includegraphics[width = 0.8 \textwidth]{ 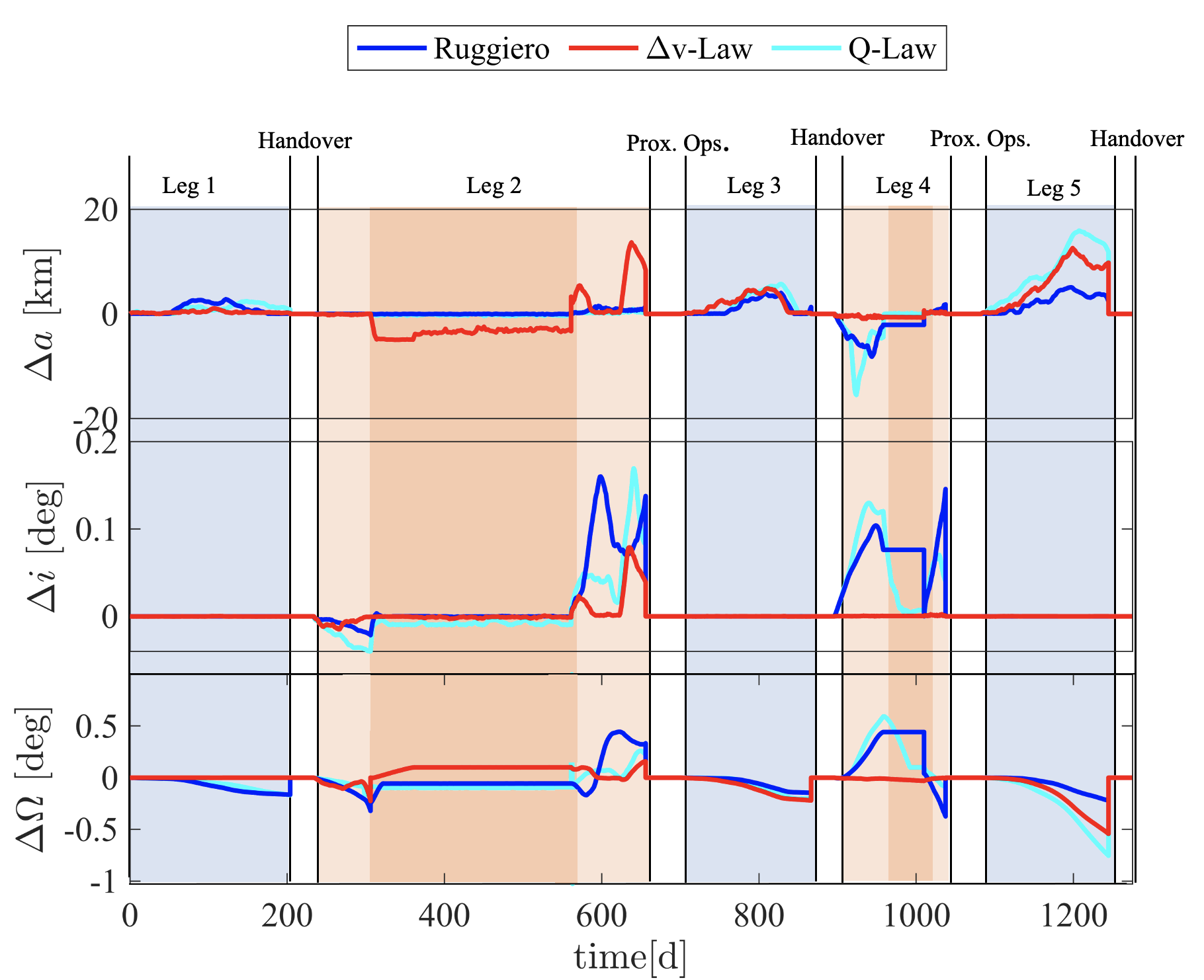}
    \caption{Time optimal tour guidance errors}
    \label{error_topt}
\end{figure}

\begin{table}[hbt!]
\centering
\caption{Time optimal tour: Error per leg obtained for each of the guidance schemes used}
\label{error_topt_table}
\begin{tabular}{c|ccc|ccc}
Leg no. & \multicolumn{3}{c|}{Forward Propagated PMDT}               & \multicolumn{3}{c}{Ruggeiro guidance}                      \\
        & $\Delta a$ (km) & $\Delta i$ (deg) & $\Delta \Omega$ (deg) & $\Delta a$ (km) & $\Delta i$ (deg) & $\Delta \Omega$ (deg) \\ \hline
1       & 2.9424          & 0.000044         & 0.3086                & 0.0847          & 0.000004         & 0.1641                \\
2       & 3.7791          & 1.315875         & 4.6885                & 0.0017          & 0.018552         & 0.2086                \\
3       & 6.7373          & 0.000022         & 0.3806                & 1.3395          & 0.000037         & 0.1480                \\
4       & 0.1680          & 0.016903         & 0.0190                & 0.0135          & 0.000815         & 0.0573                \\
5       & 10.0930         & 0.000011         & 0.6387                & 2.2808          & 0.000003         & 0.2212                \\ \hline
Total   & 23.7199         & 1.333            & 6.0354                & 3.7201          & 0.019            & 0.7991                \\ \hline \bottomrule
Leg no. & \multicolumn{3}{c|}{$\Delta v$-Law Guidance}               & \multicolumn{3}{c}{Q-Law Guidance}                         \\
        & $\Delta a$ (km) & $\Delta i$ (deg) & $\Delta \Omega$ (deg) & $\Delta a$ (km) & $\Delta i$ (deg) & $\Delta \Omega$ (deg) \\ \hline
1       & 0.1778          & 0.000107         & 0.0000                & 0.8002          & 0.000006         & 0.1672                \\
2       & 8.4013          & 0.039387         & 0.1574                & 0.3874          & 0.084089         & 0.2247                \\
3       & 0.1220          & 0.000006         & 0.2208                & 0.9388          & 0.000015         & 0.2020                \\
4       & 0.1556          & 0.000676         & 0.0076                & 0.0195          & 0.002687         & 0.0007                \\
5       & 9.7757          & 0.000067         & 0.5419                & 11.8298         & 0.000017         & 0.0755                \\ \hline
Total   & 18.6324         & 0.040            & 0.9276                & 13.9756         & 0.087            & 0.6701                \\ \hline \bottomrule
\end{tabular}
\end{table}

\begin{table}[hbt!]
\centering 
\caption{Time optimal tour: $\Delta$v, TOF and total error comparison for the guidance strategies implemented}
\label{t32}
\begin{tabular}{cccccc}
\toprule
Case                    & TOF   (d) & Fuel   (kg) & $\Delta a$ (km) & $\Delta i$ (deg) & $\Delta \Omega$ (deg) \\ \hline
Extended Edelbaum Ref.  & 1274.5    & 166.18      &23.720 & 	1.333 &	6.035 \\
Ruggiero guidance       & 1274.4    & 168.87      & 3.720           & 0.019            & 0.799                 \\
$\Delta v$-Law guidance & 1274.2    & 167.67      & 18.632          & 0.040            & 0.928                 \\
Q-Law guidance          & 1274.6    & 167.10      & 13.976          & 0.087            & 0.670                 \\ \bottomrule
\end{tabular}
\end{table}

\section{Conclusions} \label{sec:C}
The paper details the design and guidance of a multi-ADR mission. Firstly, the proposed mission architecture is discussed in detail. Then, a  preliminary mission design tool (PMDT) that considers the effect of drag, eclipses, duty ratio, and $J_2$ perturbations is developed to analyze the multi-ADR mission.  Guidance algorithms are introduced to assess the PMDT accuracy and to propose a method to track the reference trajectories. Example time and fuel optimization cases are provided for a three-debris removal mission. 

The example optimizations show that the simplified models adopted in the PMDT produce good estimates of the time of flight and propellant usage for a complex ADR mission. Furthermore, the three guidance schemes can track with good accuracy the reference trajectories even when osculating dynamics and realistic operational constraints are accounted for. 
 Hence, it is shown the method developed can optimize multi-ADR missions with a good degree of accuracy and limited computational cost. By tracking the reference trajectories computed with the PMDT it is shown that guidance laws can effectively exploit $J_2$ perturbation to reduce the propellant cost, a feature that was not considered in their original formulation. 
 
 The PMDT performance is expected to degrade with the inclusion of additional perturbations and errors such as thrust execution and orbit determination. However, the guidance laws shall provide similar accuracy in higher fidelity dynamics due to the feedback mechanisms present. Simulations done with high-fidelity dynamics are beyond the scope of this work and shall be explored in the future.  

\section*{Acknowledgments}
\label{acknowledgements}
This project was partially supported by the Ministry of Business, Innovation and Employment (MBIE) study: Astroscale/ Rocket Lab/ Te Pūnaha Ātea-Space Institute Active Debris Removal Study.

\bibliographystyle{unsrt}

\end{document}